\documentclass[
final
, nomarks
]{dmtcs-episciences}

\usepackage[utf8]{inputenc}
\usepackage{subfigure}

\usepackage{amsthm}
\usepackage{amsmath}
\usepackage{units} % for \nicefrac
\usepackage{tikz}

\newtheorem{theorem}{Theorem}
\newtheorem{proposition}[theorem]{Proposition}
\newtheorem{corollary}[theorem]{Corollary}

\usepackage[numbers]{natbib}

\author{Vincent Vatter\affiliationmark{1}\thanks{The author was partially supported by the National Science Foundation under Grant Number DMS-1301692.}}
\title{An Erd\H{o}s--Hajnal analogue for permutation classes}
\affiliation{
  Department of Mathematics, University of Florida, Gainesville, Florida, USA}
\received{2015-11-4}
\accepted{2016-2-8}
\begin{document}
\publicationdetails{18}{2016}{2}{4}{1328}
\maketitle
\begin{abstract}
Let $\mathcal{C}$ be a permutation class that does not contain all layered permutations or all colayered permutations. We prove that every permutation in $\mathcal{C}$ contains a monotone subsequence of linear size.
\end{abstract}

%
%
%
%

% Commands:
\newcommand{\Av}{\operatorname{Av}}
\newcommand{\Grid}{\operatorname{Grid}}
\newcommand{\C}{\mathcal{C}}
\newcommand{\D}{\mathcal{D}}
\newcommand{\st}{\::\:}
% Points of absolute size (to be used when drawing matchings and permutations):
\newcommand\absdot[2]{
	% Make a dot of fixed absolute size.
	\node at #1 {\normalsize $\bullet$};
	\node at #1 [below] {$#2$};
}
% Plotting permutations. Recommended scale: 0.25 (and possibly up to 0.5).
\newcommand{\plotperm}[1]{
	\foreach \j [count=\i] in {#1} {
		\absdot{(\i,\j)}{};
	};
}
\newcommand{\plotpermbox}[4]{
	\draw [thick, line cap=round]
		({#1-0.5}, {#2-0.5}) rectangle ({#3+0.5}, {#4+0.5});
}

%
%
%
%

% Get a symbol footnote:
\setcounter{footnote}{1}
\renewcommand\thefootnote{\fnsymbol{footnote}}

Ramsey's Theorem tells us that every graph on $n$ vertices has a clique or independent set of size at least a constant times $\log n$. In 1989, Erd\H{o}s and Hajnal conjectured that one can do much better when avoiding a fixed induced subgraph:

\newtheorem*{erdos-hajnal-conjecture}{The Erd\H{o}s--Hajnal Conjecture~\cite{erdos:ramsey-type-the:}}
\begin{erdos-hajnal-conjecture}
\emph{For every graph $H$ there is a constant $\delta>0$ such that every graph on $n$ vertices which does not contain $H$ as an induced subgraph has a clique or independent set of size at least $n^\delta$.}
\end{erdos-hajnal-conjecture}

We refer to Chudnovsky's survey~\cite{chudnovsky:the-erdos--hajn:} for an account of the progress that has been made toward the Erd\H{o}s--Hajnal Conjecture. Here we investigate a version of this conjecture for permutations. In this context, the analogue of the induced subgraph order is the permutation pattern order and the analogues of cliques and independent sets are decreasing and increasing subsequences, respectively. Thus the naive translation of the Erd\H{o}s--Hajnal Conjecture to permutations is well known to hold with $\delta=\nicefrac{1}{2}$:

\newtheorem*{erdos-szekeres-theorem}{The Erd\H{o}s--Szekeres Theorem~\cite{erdos:a-combinatorial:}}
\begin{erdos-szekeres-theorem}
\emph{Every permutation of length at least $(k-1)(\ell-1)+1$ contains either a decreasing subsequence of length $k$ or an an increasing subsequence of length $\ell$.}
\end{erdos-szekeres-theorem}

We instead ask for a stronger conclusion. A \emph{permutation class} is a downset of permutations under the containment order. We show (Theorem~\ref{thm-EH-main}) that unless a permutation class contains all layered permutations or all colayered permutations, its members contain monotone subsequences of linear size.

The Erd\H{o}s--Szekeres Theorem is best possible in that there are permutations of length $(k-1)(\ell-1)$ which have neither decreasing subsequences of length $k$ nor increasing subsequences of length $\ell$. Focusing on the case $k=\ell$, we call a permutation of length $k^2$ \emph{Erd\H{o}s--Szekeres extremal} if it does not contain a monotone subsequence of length $k+1$. If follows from the Robinson--Schensted correspondence and the Hook Length Formula that the number of Erd\H{o}s--Szekeres extremal permutations of length $k^2$ is
\[
	\left(\frac%
	{(k^2)!}%
	%{\left(1\cdot (2k-1)\right)\ \left(2\cdot (2k-2)\right)^2\ \left(3\cdot(2k-3)\right)^3\ \cdots\ \left((k-1)\cdot (k+1)\right)^{k-1}\ k^k}%
	{1^1\cdot 2^2\cdot 3^3\cdots k^k\cdot (k+1)^{k-1}\cdot (k+2)^{k-2}\cdots (2k-1)^1}%
	\right)^2,
\]
sequence A079402 in the OEIS~\cite{sloane:the-on-line-enc:}. One of the earliest references to this formula is Stanley's \emph{Monthly} solution~\cite{stanley:solution-to-pro:5641} from 1969. There has also been some recent interest in these permutations. In particular, Romik~\cite{romik:permutations-wi:} computed the limit shape of Erd\H{o}s--Szekeres extremal permutations. Additional results in this direction have been obtained by Kim~\cite{kim:on-increasing-s:}, Pittel and Romik~\cite{pittel:limit-shapes-fo:}, and Su~\cite{su:on-increasing-s:}.

\begin{figure}
\begin{center}
    \begin{tikzpicture}[scale=.25]
		\plotperm{3,2,1,6,5,4,9,8,7};
		\plotpermbox{1}{1}{3}{3};
		\plotpermbox{4}{4}{6}{6};
		\plotpermbox{7}{7}{9}{9};
    \end{tikzpicture} 
\quad\quad\quad\quad
    \begin{tikzpicture}[scale=.25]
		\plotperm{7,8,9,4,5,6,1,2,3};
		\plotpermbox{1}{7}{3}{9};
		\plotpermbox{4}{4}{6}{6};
		\plotpermbox{7}{1}{9}{3};
    \end{tikzpicture} 
\end{center}
\caption{The layered permutation $321\ 654\ 987=\oplus^3 (321)$ and the colayered permutation $789\ 456\ 123=\ominus^3 (123)$.}
\label{fig-layered-EES}
\end{figure}
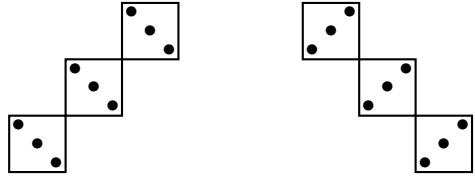

Despite the fact that there are a great many extremal Erd\H{o}s--Szekeres permutations, we need to define only two families, the \emph{layered} and \emph{colayered} extremal Erd\H{o}s--Szekeres permutations, of which examples are plotted in Figure~\ref{fig-layered-EES}. It is convenient to use the \emph{sum} and \emph{skew sum} operations in defining these permutations, which can be defined pictorially as below.
\begin{center}
	$\pi\oplus\sigma=$
	\begin{tikzpicture}[scale=0.5, baseline=(current bounding box.center)]
		\draw[thick, line cap=round] (0,0) rectangle (1,1);
		\draw[thick, line cap=round] (1,1) rectangle (2,2);
		\node at (0.5,0.5) {$\pi$};
		\node at (1.5,1.5) {$\sigma$};
	\end{tikzpicture}
\quad\quad\quad\quad
	$\pi\ominus\sigma=$
	\begin{tikzpicture}[scale=0.5, baseline=(current bounding box.center)]
		\draw[thick, line cap=round] (0,1) rectangle (1,2);
		\draw[thick, line cap=round] (1,0) rectangle (2,1);
		\node at (0.5,1.5) {$\pi$};
		\node at (1.5,0.5) {$\sigma$};
	\end{tikzpicture}
\end{center}
We use the symbols $\oplus^k$ and $\ominus^k$ to denote repeated sums and skew sums, respectively, so the extremal Erd\H{o}s--Szekeres layered permutation of length $k^2$ is
\[
	\oplus^k (k\cdots 21)
	=
	\underbrace{(k\cdots 21)\oplus\cdots\oplus (k\cdots 21)}_{\text{\footnotesize $k$ summands}}.
\]
We show that these two families of extremal permutations are the only obstructions preventing a class from having a linear bound on monotone subsequences:

\begin{theorem}
\label{thm-EH-main}
Let $\C$ be a permutation class that does not contain the class of layered permutations or the class of colayered permutations. There is a constant $c>0$ such that every permutation of length $n$ in $\C$ contains a monotone subsequence of length at least $cn$.
\end{theorem}

One special case of Theorem~\ref{thm-EH-main} is quite easy. By the folklore result that every permutation avoiding $k\cdots 21$ can be written as the union of $k-1$ increasing sequences (which can also be viewed as a consequence of Greene's Theorem~\cite{greene:an-extension-of:}), we obtain the following.

\begin{proposition}
\label{prop-EH-mono}
Let $\C$ be a permutation class that does not contain $k\cdots 21$ (or, by symmetry, $12\cdots k$). Then every permutation of length $n$ in $\C$ contains a monotone subsequence of length at least $n/(k-1)$.	
\end{proposition}

To establish Theorem~\ref{thm-EH-main} we first recall some results on transversals of (nonhomogeneous\footnote{There is also a homogeneous case, for which we refer to Kaiser~\cite{kaiser:transversals-of:}.}) $d$-intervals (although we only appeal to the case $d=2$ in what follows). Let $L_1$, $\dots$, $L_d$ be disjoint copies of the real line. A $d$-interval $I$ is a nonempty subset of $L_1\cup\cdots\cup L_d$ such that each \emph{component} $I\cap L_i$ is either empty or a closed interval. The \emph{packing number} of a set $\mathfrak{I}$ of $d$-intervals is the maximum number of pairwise disjoint elements of $\mathfrak{I}$, and this quantity is denoted $\nu(\mathfrak{I})$. A subset of $L_1\cup\cdots \cup L_d$ which intersects every member of $\mathfrak{I}$ is called a \emph{transversal} of $\mathfrak{I}$ and the minimum size of a transversal of $\mathfrak{I}$ is its \emph{transversal number}, $\tau(\mathfrak{I})$.

\begin{figure}
\begin{footnotesize}
\begin{center}
    \begin{tikzpicture}[scale=.35]
		% L_x:
    	\draw [thick, dotted] (4,0)--(4,4);
    	\absdot{(4,0)}{};
		\draw [<->, very thick] (0,0)--(11,0);
		\draw [thick] (1,1)--(4,1);
		\draw [thick] (3,2)--(6,2);
		\draw [thick] (7,3)--(10,3);
		\node [left] at (0,0) {$L_x$};
		\node [above] at (2,1) {$I_1$};
		\node [above] at (5,2) {$I_2$};
		\node [above] at (8.5,3) {$I_3$};
		% L_y:
    	\draw [thick, dotted] (5,6)--(5,10);
    	\absdot{(5,6)}{};
		\draw [<->, very thick] (0,6)--(11,6);
		\draw [thick] (1,7)--(4,7);
		\draw [thick] (5,8)--(8,8);
		\draw [thick] (3,9)--(6,9);
		\node [left] at (0,6) {$L_y$};
		\node [above] at (2.5,7) {$I_1$};
		\node [above] at (7,8) {$I_2$};
		\node [above] at (4,9) {$I_3$};
    \end{tikzpicture} 
\quad\quad\quad\quad
    \begin{tikzpicture}[scale=.35]
    	\draw [thick, dotted] (3.75,0)--(3.75,10);
    	\draw [thick, dotted] (4.25,0)--(4.25,10);
    	\draw [thick, dotted] (-1,4.75)--(11,4.75);
    	\draw [thick, dotted] (-1,5.25)--(11,5.25);
		\draw [<->, very thick] (0,0)--(11,0);
		\draw [<->, very thick] (-1,1)--(-1,10);
		\draw [thick] (1,2) rectangle (4,5);
		\draw [thick] (3,6) rectangle (6,9);
		\draw [thick] (7,4) rectangle (10,7);
		\node at (2.5,3.5) {$1$};
		\node at (5.25,7.5) {$2$};
		\node at (8.5,6.25) {$3$};
		\node [above] at (-1,10) {$L_y$};
		\node [right] at (11,0) {$L_x$};
    \end{tikzpicture} 
\end{center}
\caption{A set of $2$-intervals (left) arising from a set of axis-parallel rectangles (right).}
\label{fig-2-interval}
\end{footnotesize}
\end{figure}
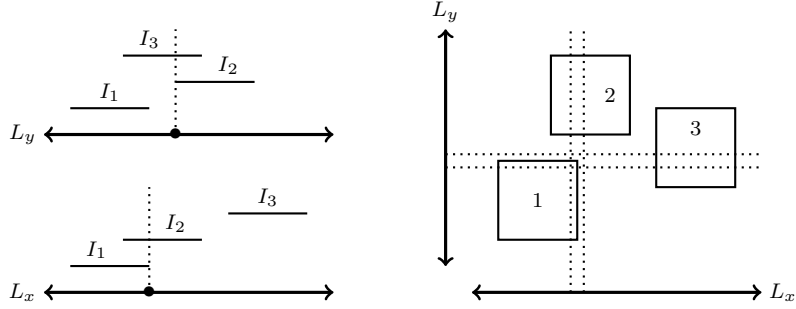

An example of a set of $2$-intervals is shown on the left of Figure~\ref{fig-2-interval}. In this example, the two copies of the real line are $L_x$ and $L_y$ and the $2$-intervals are $\mathfrak{I}=\{I_1,I_2,I_3\}$. We have $\nu(\mathfrak{I})=1$ because the $2$-intervals are pairwise intersecting, while $\tau(\mathfrak{I})=2$ because $2$ points of $L_x\cup L_y$ are required to intersect each of $I_1$, $I_2$, and $I_3$.

It is trivial to see that $\tau(\mathfrak{I})\ge \nu(\mathfrak{I})$. In the other direction, building on and improving work of Gy\'arf\'as and Lehel~\cite{gyarfas:a-helly-type-pr:} and K\'arolyi and Tardos~\cite{karolyi:on-point-covers:}, Tardos established the following.

\begin{theorem}[Tardos~\cite{tardos:transversals-of:}]
\label{thm-2int-transversal}
For every system $\mathfrak{I}$ of $2$-intervals, $\tau(\mathfrak{I})\le 2\nu(\mathfrak{I})$.
\end{theorem}

Note that the bound in Theorem~\ref{thm-2int-transversal} is best possible given the example from Figure~\ref{fig-2-interval}.

Henceforth we adopt a geometric viewpoint of permutations where a permutation $\pi$ is identified with its \emph{plot}, the set of points $\{(i,\pi(i))\}$ in the plane. Given two intervals $X$ and $Y$ of real numbers, we denote by $X\times Y$ the axis-parallel rectangle consisting of the points $\{(x,y)\st x\in X, y\in Y\}$. We further define $\pi(X\times Y)$ to be the subpermutation of $\pi$ formed by the entries of $\pi$ whose indices lie in $X$ and values lie in $Y$.

Given a set of closed axis-parallel rectangles $\mathfrak{R}$ we can naturally associate to it a set $\mathfrak{I}$ of $2$-intervals in the following manner. Let $L_x$ and $L_y$ be disjoint copies of the real line. For each $R\in\mathfrak{R}$ we create a $2$-interval $I_R$ such that $I_R\cap L_x$ is the projection of $R$ onto the $x$-axis and $I_R\cap L_y$ is the projection of $R$ onto the $y$-axis. The set of three rectangles shown on the right of Figure~\ref{fig-2-interval} corresponds in this manner to the set of $2$-intervals shown on the left of the figure.

We say that two rectangles are {\it independent\/} if both their $x$- and $y$-axis projections are disjoint. A set of rectangles is said to be independent if they are pairwise independent, and we define the \emph{independence number} of a set of axis-parallel rectangles to be the size of the largest independent set of rectangles it contains.

It follows readily that if $\mathfrak{R}$ is a set of closed axis-parallel rectangles and $\mathfrak{I}$ is the corresponding set of $2$-intervals then the independence number of $\mathfrak{R}$ is equal to the packing number of $\mathfrak{I}$. Therefore Theorem~\ref{thm-2int-transversal} shows that $\mathfrak{I}$ has a transversal of size at most $2\nu(\mathfrak{I})$. That is, there is a set of at most $2\nu(\mathfrak{I})$ points of $L_x\cup L_y$ which intersect every $2$-interval in $\mathfrak{I}$. We want to convert this transversal of $\mathfrak{I}$ into a set of vertical and horizontal lines that \emph{slice} (i.e., intersect the interior) of every rectangle in $\mathfrak{R}$, but we must be careful because the transversal is allowed to contain endpoints of $2$-intervals. Thus given a transversal $T\subseteq L_x\cup L_y$ we associate to it a set of vertical and horizontal lines $\mathcal{L}_T$ in the following way: if $t\in T$ is from $L_x$ we associate the two vertical lines $x=t\pm\epsilon$ while if $t\in T$ is from $L_y$ we associate the two horizontal lines $y=t\pm\epsilon$ where $\epsilon>0$ is a small constant (depending on $\mathfrak{R}$). These are the dotted lines shown in Figure~\ref{fig-2-interval}. This construction yields the following result.

\begin{corollary}
\label{cor-2int-slice}
If a finite collection $\mathfrak{R}$ of axis-parallel rectangles has independence number at most $m$ then there is a collection of $4m$ vertical and horizontal lines which slice every rectangle in $\mathfrak{R}$.
\end{corollary}

We need one more definition. An {\it increasing sequence\/} of rectangles is a sequence $R_1,\dots,R_m$ of independent rectangles such that $R_{i+1}$ lies above and to the right of $R_i$ for all $i$. {\it Decreasing sequences\/} of rectangles are defined analogously. We are now ready to prove Theorem~\ref{thm-EH-main}.

\begin{proof}[of Theorem~\ref{thm-EH-main}]
Let $\C$ be a permutation class which does not contain the class of layered permutations or the class of colayered permutations. Thus there is an integer $k$ such that $\C$ contains neither $\oplus^k (k\cdots 21)$ nor $\ominus^k (12\cdots k)$. Now choose a permutation $\pi\in\C$ and let $\mathfrak{R}_\pi$ denote the collection of all axis-parallel rectangles $R$ such that $\pi(R)$ contains both $k\cdots 21$ and $12\cdots k$. If $\mathfrak{R}_\pi$ were to contain an independent set of $(k-1)^2+1$ rectangles then the Erd\H{o}s--Szekeres Theorem would imply that $\mathfrak{R}_\pi$ contained an increasing or decreasing sequence of $k$ rectangles, but that would imply that $\pi$ (and hence also $\C$) contained $\oplus^k (k\cdots 21)$ or $\ominus^k (12\cdots k)$, a contradiction. Thus the independence number of $\mathfrak{R}_\pi$ is at most $(k-1)^2$.

Corollary~\ref{cor-2int-slice} now shows that there is a collection of at most $4(k-1)^2$ vertical and horizontal lines which slice every rectangle in $\mathfrak{R}_\pi$. These lines divide the plot of $\pi$ into a $t\times u$ grid of unsliced regions for integers $t$ and $u$ satisfying $t+u\le 4(k-1)^2+2$. Moreover, each unsliced region avoids either $k\cdots 21$ or $12\cdots k$ and thus can be expressed as the union of $k-1$ increasing or decreasing permutations by Proposition~\ref{prop-EH-mono}. All of the entries of $\pi$ in regions avoiding $k\cdots 21$ can be expressed as the union of at most $(t+u-1)(k-1)\le 4k^3$ increasing permutations while the entries lying in regions avoiding $12\cdots k$ can be expressed as the union of the same number of decreasing permutations. Therefore $\pi$ itself can be expressed as the union of at most $8k^3$ monotone permutations and thus $\pi$ must contain a monotone subsequence of length at least $n/8k^3$, proving the theorem.
\end{proof}

In the language of generalized grid classes introduced in Vatter~\cite{vatter:small-permutati:}, the proof of Theorem~\ref{thm-EH-main} shows that all permutation classes that do not contain the layered or colayered permutations are $\Av(k\cdots 21)\cup\Av(12\cdots k)$-griddable. In fact, this follows from \cite[Theorem 3.1]{vatter:small-permutati:}, though the bounds given in the proof above are better than that theorem gives.

Theorem~\ref{thm-EH-main} can also be interpreted in the context of splittability, which has recently been studied in general by Jel{\'{\i}}nek and Valtr~\cite{jelinek:splittings-and-:}. Given two permutation classes $\C$ and $\D$, their \emph{merge}, $\C\odot\D$, consists of those permutations whose entries can be partitioned into two subsequences, one order isomorphic to a permutation in $\C$ and the other order isomorphic to a permutation in $\D$. Thus Theorem~\ref{thm-EH-main} shows that every class which does not contain the layered or colayered classes can be expressed as the merge of a finite number of monotone classes (i.e., classes of the form $\Av(21)$ or $\Av(12)$). Note that the converse---if a class does contain all layered or colayered permutations then it is not the merge of a finite number of monotone classes---holds trivially. The proof of Theorem~\ref{thm-EH-main} also gives a bound on how many monotone classes are required:

\begin{corollary}
\label{cor-merge-mono}
If a permutation class avoids both $\oplus^k (k\cdots 21)$ and $\ominus^k (12\cdots k)$ then it can be expressed as the merge of at most $8k^3$ monotone classes.
\end{corollary}

%\nocite{*}
%\bibliographystyle{abbrvnat}
% use the following instead if you encounter problems 
%\bibliographystyle{alpha}
\bibliographystyle{abbrv}
\bibliography{../../../refs}
\label{sec:biblio}

\end{document}